\documentclass[12pt,twoside,reqno]{amsart}
\usepackage{pifont}
\usepackage{amsmath}
\usepackage{amsthm}
\usepackage{amsfonts}
\usepackage{amssymb}
\usepackage{latexsym}
\usepackage{amstext}
\usepackage{array}
\usepackage{graphicx}

\date{}
\allowdisplaybreaks[4] \footskip=15pt
\renewcommand{\uppercasenonmath}[1]{}

\newtheorem{thm}[subsection]{Theorem}
\newtheorem{cor}[subsection]{Corollary }
\newtheorem{Def}[subsection]{Definition}
\newtheorem{lem}[subsection]{Lemma}
\newtheorem{remark}[subsection]{Remark}

\newtheorem{prop}[subsection]{Proposition}
\newtheorem{exm}[subsection]{Example}
  
\newcommand{\bthm}{\begin{thm} }
\newcommand{\ethm}{\end{thm} }
\newcommand{\bpro}{\begin{prop}}
\newcommand{\epro}{\end{prop}}
\newcommand{\bdf}{\begin{Def}}
\newcommand{\edf}{\end{Def}}
\newcommand{\bexm}{\begin{exm}}
\newcommand{\eexm}{\end{exm}}
\newcommand{\blem}{\begin{lem}}
\newcommand{\elem}{\end{lem}}
\newcommand{\bpf}{\begin{proof}}
\newcommand{\epf}{\end{proof}}
\newcommand{\bcor}{\begin{cor}}
\newcommand{\ecor}{\end{cor}}
\newcommand{\ba}{\begin{array}}
\newcommand{\ea}{\end{array}}
\newcommand{\bea}{\begin{eqnarray}}
\newcommand{\eea}{\end{eqnarray}}

\newcommand{\brem}{\begin{remark}}
\newcommand{\erem}{\end{remark}}

\begin{document}

\begin{center}
{\large  \bf Lehmer's totient problem over $\mathbb{F}_q[x]$}
\footnote {Supported by NSFC (Nos. 11271177, 11471154, 11571163).}\\

 \vskip 0.8cm
 {\small   Qingzhong Ji \footnote{Corresponding author.\\ \indent E-mail addresses: qingzhji@nju.edu.cn (Q. Ji),\; hrqin@nju.edu.cn (H. Qin)} \ \ and\ \     Hourong Qin}\\
{\small Department of Mathematics, Nanjing University, Nanjing
210093, P.R.China}

\vskip 3mm
\end{center}

{\bf Abstract:}  {\small In this paper, we consider the function field analogue of the Lehmer's totient problem. Let $p(x)\in\mathbb{F}_q[x]$ and $\varphi(q,p(x))$ be  the  Euler's totient function of $p(x)$ over  $\mathbb{F}_q[x],$ where  $\mathbb{F}_q$ is a finite field with $q$ elements. We prove that $\varphi(q,p(x))|(q^{{\rm deg}(p(x))}-1)$ if and only if (i) $p(x)$ is irreducible;  or (ii) $q=3, \; p(x)$ is the product of any $2$  non-associate irreducibes of degree $1;$ or (iii)  $q=2,\; p(x)$ is the product of all irreducibles of degree $1,$ all irreducibles of degree $1$ and $2,$ and the product of any $3$ irreducibles one each of degree $1, 2$ and $3$.}



\section{\bf Introduction}\label{1}

Throughout this paper, let $\Bbb Q$, $\Bbb Z$ and $\Bbb N$ denote
the field of rational numbers, the ring of rational integers and the
set of nonnegative integers, respectively. Let ${\Bbb N}^*={\Bbb N}\setminus \{
0\}$. As usual, let ${\rm ord}_p$ denote the normalized $p$-adic valuation
of $\Bbb Q_p.$

{\bf Lehmer's totient problem }   Let $\varphi$ be  the Euler's totient function.  In \cite{Le}, Lehmer discussed the equation \begin{eqnarray}\label{10001}k\varphi(n)=n-1,\end{eqnarray} where $k$ is an integer. In his pioneering  paper \cite{Le}, Lehmer showed that if $n$ is a solution of (\ref{10001}), then $n$ is a prime or the product of seven or more distinct primes.  One is tempted to believe that an integer $n$ is a prime if and only if $\varphi(n)$ divides $n-1.$ This problem  has not been solved to this day. But some progress has been made in this direction. In the literature, some authors call these composite numbers $n$ satisfying  equation (\ref{10001}) the Lehmer numbers. Lehmer's totient problem is to determine the set of  Lehmer numbers. To the best of our knowledge, the current best result is due to Richard G. E. Pinch(see\cite{Pin}), that the number of prime factors of a Lehmer number $n$ must be at least 15 and there is no Lehmer number less than $10^{30}.$   For further results on this topic we refer the reader to  (\cite{BGW}, \cite{BF}, \cite{GO}, \cite{LP}, \cite{Po}).

J. Schettler \cite{Sc} generalizes the divisibilty condition $\varphi(n)|(n-1),$ constructs reasonable notion of Lehmer numbers and Carmichael numbers in a PID and gets some interesting results.  Let $R$ be a PID with the property:  $R/(r)$ is finite whenever $0\not= r\in R.$  Denote the sets of units, primes and (non-zero) zero divisors, in $R,$ by $U(R),\;P(R)$ and $Z(R),$ respectively; additionally, define
\begin{eqnarray}\label{10002}L_R:=\{r\in R\setminus (\{0\}\cup U(R)\cup P(R)):\; |U(R/(r))|\;|\;|Z(R/(r))|\}.\end{eqnarray}
Note that when $R=\mathbb{Z},$  $L_{\mathbb{Z}}$ is the set of Lehmer numbers. An element of $L_R$ is also called a Lehmer number of $R.$ Let  $\mathbb{F}_q$ is a finite field with $q$ elements. Then $\mathbb{F}_q[x]$ is a PID.  Schettler obtains some properties of elements of $L_{\mathbb{F}_q[x]}$ as  follows.

\bpro {\rm(\cite{Sc}, Theorems 5.1, 5.2, 5.3 )}  {\rm (1)}  Suppose $f(x)\in L_{\mathbb{F}_q[x]},$ $p(x)\in P(\mathbb{F}_q[x])$ and $p(x)|f(x).$  Then ${\rm deg}(p(x))|{\rm deg}(f(x)).$

{\rm (2)}     Suppose $f(x)\in L_{\mathbb{F}_q[x]}.$   Then $f(x)$ has at least $[{\rm log}_2(q+1)]$ distinct prime factors.

 {\rm (3)}    There exists a PID $R$ such that $L_R\not=\emptyset.$  {\rm(}{\sl E.g.}, $f(x)=x(x+1)\in L_{\mathbb{Z}/{2\mathbb{Z}}.}{\rm)}$
\epro

\vskip 2mm

Our work is inspired by  above proposition,  in this paper, our goal is to determine the set $L_{\mathbb{F}_q[x]}.$

\vskip 2mm

{\bf Euler's totient function over $\mathbb{F}_q[x]$.}   Let $f(x)\in\mathbb{F}_q[x]$ with $m={\rm deg}(f(x))\geq 1.$  Put
 $$\Phi(f(x))=\{g(x)\in\mathbb{F}_q[x]\;|\;{\rm deg}(g(x))\leq m-1,(f(x),g(x))=1\}.$$ The Euler's totient function $\varphi(q,f(x))$ of $f(x)$ is defined as follows:
$$\varphi(q,f(x))=\sharp \Phi(f(x)).$$
If $f(x)\in \mathbb{F}_q[x]$ is irreducible, then   $\varphi(q,f(x))=q^{{\rm deg}(f(x))}-1.$ It is easy to see that  the functions $\varphi(q,f(x))$ and $\varphi(n)$ have the following similar  properties:

\bpro \label{001}Let $f(x)=p_1(x)^{r_1}\cdots p_k(x)^{r_k}\in\mathbb{F}_q[x]$ of degree $n\geq 1,$ where $p_1(x),\ldots,p_k(x)\in P(\mathbb{F}_q[x])$ are non-associate,  ${\rm deg}(p_i(x))=n_i$ and $r_i\geq 1,$  $1\leq i\leq k.$ Then we have

{\rm (1)}  $\varphi(q,f(x))=q^n\prod\limits_{i=1}^k(1-\frac{1}{q^{n_i}});$

{\rm (2)} If $g(x)\in \mathbb{F}_q[x]$ and $(f(x),g(x))=1,$ then $g(x)^{\varphi(q,f(x))}\equiv 1\pmod{f(x)};$

{\rm (3)} If $\varphi(q,f(x))|(q^n-1),$ then $r_i=1,$ for all $1\leq i\leq k.$\epro

Hence it is  natural to consider the Lehmer's totient problem over $\mathbb{F}_q[x]$:
\begin{center} Determine $f(x)\in\mathbb{F}_q[x]$ such that $\varphi(q,f(x))|(q^{{\rm deg}(f(x))}-1).$ \end{center}
Set $$\mathcal{L}_{\mathbb{F}_q}=\{f(x)\in \mathbb{F}_q[x]\setminus\{0\}\;|\;{\rm deg}(f(x))\geq 1,\;\varphi(q,f(x))|(q^{{\rm deg}(f(x))}-1)\}.$$
By the definition (\ref{10002}), it is easy to see  that $$L_{\mathbb{F}_q[x]}=\{f(x)\in \mathbb{F}_q[x]\setminus\{0\}\;|\;f(x)\;{\rm is\;reducible},\;\varphi(q,f(x))|(q^{{\rm deg}(f(x))}-1)\}.$$  Hence  $\mathcal{L}_{\mathbb{F}_q}=P(\mathbb{F}_q[x])\cup L_{\mathbb{F}_q[x]}.$

  For $q=2,3,$ Lv Hengfei \cite{Lv}  gave some  polynomials $f(x)\in L_{\mathbb{F}_q[x]}$ as follows:

(1) $q=2,$ $f(x)=x(x+1)(x^2+x+1),$ then $\varphi(2,f(x))=3,$ hence $\varphi(2,f(x))|(2^4-1).$

(2) $q=3,$ $f(x)=x(x+1),$ then $\varphi(3,f(x))=4,$ hence $\varphi(3,f(x))|(3^2-1).$

\vskip 2mm

In this paper, we   give the necessary and sufficient conditions for $f(x)\in L_{\mathbb{F}_q[x]}$ as follows.
\vskip 2mm
{\bf Main Theorem} {\sl \label{01} {\rm (1)} Assume $q\geq 4.$ Then $L_{\mathbb{F}_q[x]}=\emptyset.$

{\rm (2)}  Assume $q=3.$  Then $L_{\mathbb{F}_3[x]}$ consists of the products of any $2$  non-associate irreducibes of degree $1,$ {\sl i.e.},    $$L_{\mathbb{F}_3[x]}=\{ax(x+1), ax(x-1), a(x+1)(x-1)\in\mathbb{F}_3[x], a=1,2\}.$$

 {\rm(3)}  Assume $q=2.$   Then  $L_{\mathbb{F}_2[x]}$ consists of  the products of all irreducibles of degree 1, the products of all irreducibles of degree 1 and 2, and the products of any 3 irreducibles one each of degree 1, 2, and 3, {\sl i.e.},    $$\begin{array}{ll}L_{\mathbb{F}_2[x]}=&\{x(x+1), x(x+1)(x^2+x+1), x(x^2+x+1)(x^3+x+1),\\
 &\;\; (x+1)(x^2+x+1)(x^3+x+1), x(x^2+x+1)(x^3+x^2+1),\\&\;\; (x+1)(x^2+x+1)(x^3+x^2+1)\in \mathbb{F}_2[x]\}.\end{array}$$}

\vskip 2mm

The proof is essentially to give the necessary and sufficient conditions for $\varphi(q,f(x))|(q^{{\rm deg}(f(x))}-1)$  which will be divided into two cases $q\geq 3$ and  $q=2.$

{\bf Acknowledgments:} We would like to thank Prof. X. Guo and H. Lv for reading the manuscript carefully and providing valuable comments and suggestions.

\section{\bf Properties of cyclotomic polynomials}\label{2}

Let $n\in{\Bbb N}^*$ and  $\zeta_n$ be a primitive $n$-th root of unity. The polynomial
$$\Phi_n(x)=\prod\limits_{(j,n)=1}(x-\zeta_n^j)$$
is called  the $n$-th cyclotomic polynomial.
It is well-known that $\Phi_n(x)$ is an irreducible polynomial of degree $\varphi(n)$ in $\mathbb{Z}[x]$ and
\begin{eqnarray}\label{20001}x^n-1=\prod\limits_{d|n}\Phi_d(x).\end{eqnarray}
Note that the polynomial factorization in (\ref{20001}) is complete. But it does not follow that the factorization
\begin{eqnarray}\label{20002}a^n-1=\prod\limits_{d|n}\Phi_d(a),\;\;\;a\in \mathbb{Z},\end{eqnarray} is complete, since the integer $\Phi_d(a)$ may not be prime.

\bdf Suppose $a > b > 0$ are coprime integers. A prime divisor $p$ of $a^n-b^n,$ $n\geq 2,$ is called primitive if $p\nmid a^k-b^k,$ for any $k<n.$ Otherwise, it is called algebraic.\edf
It is well-known that the following Bang-Zsigmondy's Theorem provides the existence of a primitive prime factor.

\vskip 2mm

 {\bf Bang-Zsigmondy's Theorem}(\cite{Zs}) {\sl Suppose $a > b > 0$ are coprime integers. Then for any natural number $n > 1$ there is a  primitive prime divisor $p$ of  $a^n-b^n$  with the following exceptions:
$$a = 2, \;b = 1,\;{\rm  and}\;\; n = 6;\;\;{\rm or}$$
$$a + b\;\;{\rm is\;\; a\;\; power\;\; of\;\; two,\;\; and}\;\; n = 2.$$}

\vskip 2mm

It is clear that for any $n,$ and $d|n,$ that any prime $p$ dividing $\phi_d(a)$ will be an algebraic divisor of (\ref{20002}), since $p$ must divide $a^d-1$ as $\phi_d(a)$ does.  On the other hand, any primitive factor of $a^n-1$ will have to divide $\Phi_n(a).$ It is not true, however, that every prime factor of $\Phi_n(a)$ is primitive.

 \vskip 2mm

\blem \label{2002} {\rm (\cite{BLSTW}, III C1, p. lxviii)}  Let $p$ be a prime and $m\in\mathbb{N}^*$ with $(p,m)=1.$ Suppose $v\in\mathbb{N}^*$ and  $a\in\mathbb{Z}.$ Then $p|\Phi_{mp^v}(a)$ if and only if $p|\Phi_m(a).$  Furthermore,

{\rm(1)} if $p|\Phi_m(a)$ and $mp^v>2,$  then ${\rm ord}_p(\Phi_{mp^v}(a))=1;$

{\rm(2)} if  $p|\Phi_m(a)$ and $mp^v=2,$ {\sl i.e.,} $p=2, m=v=1,$   then $${\rm ord}_2(\Phi_2(a))={\rm ord}_2(a+1)\geq 1.$$
 \elem

\vskip 2mm

\blem\label{2003}
Let $p$ be a prime and $n\in\mathbb{N}^*.$ Suppose $n=p^vm$ with $v={\rm ord}_p(n).$ Then $p|\Phi_n(a)$ for some $a\in\mathbb{Z}$ if and only if $m|(p-1).$
\elem
\bpf \;\; It is obvious from Lemma \ref{2002} and (\cite{Wa}, Lemmas 2.9, 2.10).\qed
\epf

\vskip 2mm
\bcor\label{2004} Let $p$ be a prime and $a\in\mathbb{Z},\;v\in\mathbb{N}.$ Then  $p|\Phi_{p^v}(a)$  if and only if $p|(a-1).$ \ecor

\vskip 2mm

\bcor \label{2005}\; Let $m>n$ be positive integers. For any $a\in\mathbb{Z},$ we obtain that $(\Phi_n(a),\Phi_m(a))=1$ or $(\Phi_n(a),\Phi_m(a))$ is a prime. Furthermore, if $(\Phi_n(a),\Phi_m(a))=p$ is a prime,  then $m=p^vn$ for some $v\geq 1.$ \ecor

\vskip 2mm

\blem  \label{2006} Let $a,m\in\mathbb{N}^*$ and $a\geq 2.$ Then $|\Phi_m(a)|=1$ if and only if $m=1, a=2.$\elem

\bpf\;\; By the formula $\Phi_m(a)=\prod\limits_{(j,m)=1}(a-\zeta_m^j),$ we know that $|a-\zeta_m^j|>1$ for all  $a\geq 2$ and $m\geq 2,$ hence $|\Phi_m(a)|>1.$ On the other hand, $\Phi_1(x)=x-1.$  Therefore   $\Phi_m(a)=1$ if and only if $m=1,a=2.$\qed\epf
\vskip 2mm
To end this section, we  recall an estimate for $\Phi_n(a).$

\vskip 2mm
\blem {\rm (\cite{TV}, Theorem 5)}\label{TV} \label{2007} For any integers $n\geq 2$ and $a\geq 2,$ we have $$\frac{1}{2}a^{\varphi(n)}\leq \Phi_n(a)\leq 2\cdot a^{\varphi(n)}.$$ \elem

\section{\bf  Main Results}\label{3}
\vskip 0.5cm
Let the notation be the same as in \S \ref{1} and \S \ref{2}.

\vskip 2mm

\bpro \label{002} Let $a,n\in\mathbb{N}^*$  and $a\geq 3,\;n\geq 2.$ Assume $s\geq 2$ and $e_{_1},e_2,\ldots,e_s\in\mathbb{N}^*$  with $\sum\limits_{i=1}^se_i=n.$ Then $\prod\limits_{i=1}^s(a^{e_i}-1)| (a^n-1)$ if and only if

{\rm  (1)}  $a=3,$  $n=s=2,$   $e_{_1}=e_2=1,$ or

{\rm(2)} $a=3,$  $n=s=4,$  $e_1=e_2=e_3=e_4=1.$
\epro

 \bpf\;\; The sufficiency  is trivial.  It is sufficient to show the necessity.  Suppose $\prod\limits_{i=1}^s(a^{e_i}-1)| (a^n-1).$  First, we have \begin{eqnarray}\label{3001}\frac{x^n-1}{\prod\limits_{i=1}^s (x^{e_i}-1)}=\frac{\prod\limits_{d\in T}\Phi_d(x)}{ \prod\limits_{d'\in T'}\Phi_{d'}(x)}=\frac{\prod\limits_{d\in T}\Phi_d(x)}{(x-1)^{s-1}\cdot \prod\limits_{d'\in T''}\Phi_{d'}(x)}=\frac{P(x)}{Q(x)},\end{eqnarray} where $T=\{d>1\;|\;d|n,\;d\nmid e_i,\;1\leq i\leq s\},$  $P(x)=\prod\limits_{d\in T}\Phi_d(x)$ and $Q(x)= \prod\limits_{d'\in T'}\Phi_{d'}(x)$ for some index set $T',$  and  $T''=\{d'\in T'\;|\;d'\geq 2\}.$

  We have

(i) $(P(x),Q(x))=1$  and  ${\rm deg}(P(x))={\rm deg}(Q(x));$

(ii) For any $d'\in T',$ we have
  $$d'|e_i\;{\rm for\;some\;}1\leq i\leq s, \;{\rm and}\;(\Phi_{d'}(x),\Phi_d(x))=1 \;{\rm for \;all\;}d\in T;$$

  (iii) For any $d\in T$ and $d'\in T',$ we have $d\nmid d'.$

  (iv) For any $d\in T$ and $d_1',d_2'\in T'$ such that  $$(\Phi_d(a),\Phi_{d_1'}(a))\not=1\;{\rm and }\;(\Phi_d(a),\Phi_{d_2'}(a))\not=1.$$ Then $(\Phi_d(a),\Phi_{d_1'}(a))=(\Phi_d(a),\Phi_{d_2'}(a))=p$ for some prime $p$ and  $d=p^{v_1}d_1'=p^{v_2}d_2'$ for some $v_1,v_2\in\mathbb{N}^*.$ Furthermore, ${\rm ord}_p(\Phi_d(a))=1$ except $d=2, d_1'=d_2'=1.$

  The statements (i),(ii) and (iii) are obvious. We only prove (iv).  In fact, by Corollary \ref{2005}, there exist primes $p_1$ and $p_2$ such that $(\Phi_d(a),\Phi_{d_1'}(a))=p_1$ and $(\Phi_d(a),\Phi_{d_2'}(a))=p_2.$ If $p_1\not=p_2,$ then by (iii) and Corollary \ref{2005}, we have $d=p_1^{r_1}p_2^{r_2}d''$ for some $r_1,r_2,d''\in\mathbb{N}^*$ with $(p_1,p_2d'')=(p_2,p_1d'')=1.$ By  Lemma \ref{2003}, we have $p_2^{r_2}d''|(p_1-1)$ and $p_1^{r_1}d''|(p_2-1).$ This is a contradiction. Hence we obtain $(\Phi_d(a),\Phi_{d_1'}(a))=(\Phi_d(a),\Phi_{d_2'}(a))=p$ for some prime $p.$ From (iii) and Corollary \ref{2005}, we have $d=p^{v_1}d'_1=p^{v_2}d'_2$ for some $v_1,v_2\in\mathbb{N}^*.$   By Lemma \ref{2002}, we have ${\rm ord}_p(\Phi_d(a))=1$ except $d=2, d_1'=d_2'=1.$ Thus we complete the proof of (iv).

By assumption, we have
$$\frac{a^n-1}{\prod\limits_{i=1}^s (a^{e_i}-1)}=\frac{\prod\limits_{d\in T}\Phi_d(a)}{(a-1)^{s-1}\cdot \prod\limits_{d'\in T''}\Phi_{d'}(a)}=\frac{P(a)}{Q(a)}\in\mathbb{N}^*.$$

By assumption $a\geq 3,$ then either $a-1=2^r$ or there exists an odd prime $p$ such that $p^r|| (a-1)$   for some $r\in\mathbb{N}^*.$   Then $2^{r(s-1)}|P(a)$ or  $p^{r(s-1)}|P(a).$  If $a-1=2^r,$ then
$${\rm ord}_2(\Phi_2(a))={\rm ord}_2(a+1)={\rm ord}_2(2^r+2)=\left\{\begin{array}{ll}1,&r\geq 2,\\ 2,&r=1.\end{array}\right.$$

{\bf Case 1} Assume $p=2$ and $r=1,$ {\sl i.e.,} $a=3.$   Since $2|T_d$ for some $d\in T,$ $d$ is even, so is $n$ even.

(a) If $2\not\in T,$ by  Lemma \ref{2002} and  Corollary \ref{2005}, there exist positive integers $2\leq j_1<j_2<\cdots<j_{s-1}$ such that
$$2^{j_{_1}}, 2^{j_{_2}}, \ldots, 2^{j_{s-1}}\in T,\;{\rm and \;} {\rm ord}_2(\Phi_{2^{j_{_k}}}(3))=1,\;1\leq k\leq s-1.$$

(b) If $2\in T,$ then $e_1, \ldots, e_s$ are odd, hence $s$ is even.

If $s\geq 4,$  then $2, 2^2,\ldots, 2^{s-2}\in T$ and
$${\rm ord}_2(\Phi_2(3))=2, \;\;{\rm ord}_2(\Phi_{2^k}(3))=1,\;\;2\leq k\leq s-2. $$

{\bf Case 2} Assume $p$ is odd or $p=2, a-1=2^r, r\geq 2.$  By  Lemma \ref{2002} and  Corollary \ref{2005}, there exist positive integers $1\leq i_1<i_2<\cdots<i_{r(s-1)}$ such that
$$p^{i_{_1}},p^{i_{_2}},\ldots,p^{i_{r(s-1)}}\in T,\;{\rm and \;} {\rm ord}_p(\Phi_{p^{i_{_k}}}(a))=1,\;1\leq k\leq r(s-1).$$
We set  $$\Delta=\left\{\begin{array}{ll}\{2\},&{\rm if}\;p=2, a=3, 2\in T, s=2,\\
\{2, 2^2,\ldots, 2^{s-2}\},&{\rm if}\;p=2, a=3, 2\in T, s\geq 4,\\
\{2^{j_{_1}}, 2^{j_{_2}}, \ldots, 2^{j_{s-1}}\},&{\rm if}\;p=2, a=3, 2\not\in T,\\
\{p^{i_{_1}}, p^{i_{_2}}, \ldots, p^{i_{r(s-1)}}\},&{\rm if}\;p \;{\rm is\; odd\; or}\; p=2, a-1=2^r, r\geq 2.\end{array}\right.$$

If $T''\not=\emptyset,$  we define a map $f: T''\longrightarrow T$ as follows.   By Lemma \ref{2006}, for any $d'\in T'',$  we have  $|\Phi_{d'}(a)|\not=1.$  Choose a prime factor  of $\Phi_{d'}(a),$ say  $p'|\Phi_{d'}(a),$   there exists  $d=p'^vd'\in T$ for some $v\geq 1.$ Define $f(d')=d.$  By Lemma \ref{2002}, we have ${\rm ord}_{p'}(\Phi_d(a))=1.$  By (iv), the map $f$ is injective and $f(d')\not\in\Delta.$ For any $d'\in T'',$ we have $d'\geq 2$ and $p'|\Phi_{d'}(a),$ and if $p'=2,$ then $2|d'.$ Hence
$${\rm deg}(\Phi_{f(d')}(x))=\varphi(p'^vd')>\varphi(d')={\rm deg}(\Phi_{d'}(x)),\;\;d'\in T''.$$ On the other hand, we always have
$$\sum\limits_{m\in\Delta}{\rm deg}(\Phi_{m}(x))\geq s-1.$$   Hence the equality  ${\rm deg}(P(x))={\rm deg}(Q(x))$ implies that $T''=\emptyset$ and
 $$\sum\limits_{m\in\Delta}{\rm deg}(\Phi_{m}(x))= s-1\;\;{\rm and}\;\;  a-1=p^r.$$ Note that $T''=\emptyset$ implies that $$e_i|n \;{\rm and}\; (e_i, e_j)=1,\;1\leq i\not=j\leq s.$$
 It is easy to verify  that
$\sum\limits_{m\in\Delta}{\rm deg}(\Phi_{m}(x))= s-1$ if and only if (i) $a=3,$ $p=2,$ $s=2,$   $e_{_1}=e_2=1,$ or (ii) $a=3, p=2, s=4, e_1=e_2=e_3=e_4=1.$  This completes the proof.\qed\epf

\vskip 2mm

\blem\label{5001}
Let $n\in\mathbb{N}^*$  and $n\geq 2.$ Assume $s\geq 2$ and $e_{_1},e_2,\ldots,e_s\in\mathbb{N}^*$  with $\sum\limits_{i=1}^se_i=n.$ If $\prod\limits_{i=1}^s(2^{e_i}-1)| (2^n-1),$ then $e_i|n$ for all $1\leq i\leq s,$ and $(e_{_1},\ldots,e_s)=1.$\elem

\bpf\;\;The assumption  $\prod\limits_{i=1}^s(2^{e_i}-1)| (2^n-1)$ implies that
$$\frac{2^n-1}{\prod\limits_{i=1}^s (2^{e_i}-1)}=\frac{\prod\limits_{d\in T}\Phi_d(2)}{ \prod\limits_{d'\in T''}\Phi_{d'}(2)}=\frac{P(2)}{Q(2)}\in\mathbb{N}^*,$$
where the sets $T$ and $T''$ are defined by the formula (\ref{3001}). Suppose that there exists $e_{i_{_0}}$ for some $1\leq i_{_0}\leq s$ such that $e_{i_{_0}}\nmid n.$ Hence there is a prime $p$ and $r\in\mathbb{N}^*$ such that $p^r|e_{i_{_0}}$ and $p^r\nmid n.$  Thus $p^r\in T''.$  By Lemma \ref{2006}, we have $|\Phi_{p^r}(2)|\not=1.$ Let $q$ be a prime such that $q|\Phi_{p^r}(2).$ Then there exists $d\in T$ such that $q|\Phi_d(2).$  From (iii) of the proof of Proposition \ref{002} and Corollary \ref{2005},  we have $d=q^vp^r$ for some $v\in\mathbb{N}^*.$  Therefore $q^vp^r|n.$  This contradicts the fact $p^r\nmid n.$ Hence we have $e_i|n$ for all $1\leq i\leq s.$

Assume $(e_{_1},\ldots,e_s)=d>1.$  Put $a=2^d,$ $e_i=e'_id,\;1\leq i\leq s,$ $n=n'd.$  Then $a\geq 4$ and $n'=\sum\limits_{i=1}^se'_i.$ By Proposition \ref{002}, we have $\prod\limits_{i=1}^s(a^{e'_i}-1)\nmid (a^{n'}-1),$ hence  $\prod\limits_{i=1}^s(2^{e_i}-1)\nmid (2^{n}-1).$ This contradicts the assumption  $\prod\limits_{i=1}^s(2^{e_i}-1)| (2^{n}-1).$ Therefore we have $(e_{_1},\ldots,e_s)=1.$\qed\epf

\vskip 2mm

\blem\label{5002} Let
$n\in \mathbb{N}^*$ and $h(n)=\frac{\sigma(n)}{n},$ where $\sigma(n)=\sum_{d|n}d.$
Then we have $h(n)<1.28 n^{\frac{1}{4}},$ for all $n\in \mathbb{N}^*.$
 \elem
\bpf \;\;  Let $p\geq 5$ be a prime and $a\in\mathbb{N}^*.$ It is easy to see that $\frac{h(p^a)}{p^{\frac{a}{4}}}<1.$  For $p=2,\;3,$ we get
$$\frac{h(2^a)}{2^{\frac{a}{4}}}\left\{\begin{array}{ll}< 1.262,& {\rm if}\;\; a=1,\\< 1.238,&{\rm if}\;\;  a=2,\\ < 1.115,&{\rm if}\;\; a=3,\\
< 1,&{\rm if}\;\; a\geq 4.\end{array}\right.$$
and $$\frac{h(3^a)}{3^{\frac{a}{4}}}\left\{\begin{array}{ll}< 1.014,& {\rm if}\;\; a=1,\\
< 1,&{\rm if}\;\; a\geq 2.\end{array}\right.$$
Hence we have  $h(n)<1.262\times 1.014 n^{\frac{1}{4}}<1.28 n^{\frac{1}{4}},$ for all $n\in \mathbb{N}^*.$\qed\epf
\vskip 2mm
\blem \label{5003} Let  $n\in \mathbb{N}^*.$  Set
$$c(n)=\left\{\begin{array}{ll}0.59,&{\rm if}\;\;{\rm ord}_2(n)=1,\\0.70,&{\rm if}\;\;{\rm ord}_2(n)=2,\\0.84,&{\rm if}\;\; {\rm ord}_2(n)=3,\\ 1,& {\rm if}\;\;{\rm ord}_2(n)\geq 4, \;{\rm or} \;\;{\rm ord}_2(n)=0.\end{array}\right.$$
Then $\varphi(n)>c(n)n^{\frac{3}{4}},$  for any integer $n\geq 2.$\elem

\bpf\;\;If $p$ is an odd prime, then $\varphi(p^a)>p^{\frac{3a}{4}}$ for any $a\in\mathbb{N}^*.$  On the other hand, we have
$$\frac{\varphi(2^a)}{2^{\frac{3a}{4}}}\left\{\begin{array}{ll}> 0.59,&{\rm if}\;\;{\rm ord}_2(n)=1,\\>0.70,&{\rm if}\;\;{\rm ord}_2(n)=2,\\>0.84,& {\rm if}\;\;{\rm ord}_2(n)=3,\\ >1,& {\rm if}\;\;{\rm ord}_2(n)\geq 4.\end{array}\right.$$ Hence $\varphi(n)>c(n)n^{\frac{3}{4}},$  for any integer $n\geq 2.$\qed\epf

\vskip 2mm
\bpro\label{5004}
Let $n\geq s\geq 2,$    $e_1\leq e_2\leq \cdots<e_s$ be  positive integers such that $\sum\limits_{i=1}^se_i=n.$ For each $d|n,$ $d<n,$ let $u_d=\sharp \{e_i\;|\;e_i=d,\;1\leq i\leq s\}.$ Assume that
 $u_1\leq 2$ and $u_d\leq \frac{2^d-1}{d}$ for any $d\geq 2.$
  Then $\prod\limits_{i=1}^s(2^{e_i}-1)| (2^n-1)$ if and only if   ${\rm (1)}\; n=2, s=2, e_1=e_2=1;$  or ${\rm (2)}\; n=4, s=3, e_1=e_2=1,  e_3=2;$ or ${\rm (3)}\; n=6, s=3, e_1=1, e_2=2, e_3=3.$   \epro

\bpf\;\; The sufficiency  is trivial.  It is sufficient to show the necessity. Set $$R=\frac{2^n-1}{\prod\limits_{i=1}^s(2^{e_i}-1)}\in\mathbb{N}^*.$$

(1) Assume $2\leq n\leq 6.$  It is easy to show the necessity by Lemma \ref{5001}.

(2) Assume $n\geq 7.$  The primitive part $M$ of $2^n-1$ can not be reduced with the denominator, so $R\geq M.$  By Lemma  \ref{TV} ,  we have
$$R\geq M\geq \frac{\Phi_n(2)}{n}\geq \frac{2^{\varphi(n)}}{2n}.$$

On the other hand, we have $$R=\frac{2^n-1}{2^n}\prod\limits_{i=1}^s(1-2^{-e_i})^{-1}<\prod\limits_{i=1}^s(1-2^{-e_i})^{-1}.$$
By assumption, $u_1\leq 2,\;u_2\leq 1,$ hence
$$\begin{array}{ll}{\rm log}R& <2{\rm log}2+\delta(n){\rm log}\frac{4}{3}-\sum_{e_i\geq 3}{\rm log}(1-2^{-e_i})\\
&<{\rm log}4+\delta(n){\rm log}\frac{4}{3}+\sum_{e_i\geq 3} \frac{1}{2^{e_i}-1}\\
&<{\rm log}4+\delta(n){\rm log}\frac{4}{3}+\sum_{d|n,\;3\leq d<n}\frac{u_d}{2^d-1}\\
&\leq {\rm log}4+\delta(n){\rm log}\frac{4}{3}+\sum_{d|n,\;3\leq d<n}\frac{1}{d}\\
&={\rm log}4+\delta(n){\rm log}\frac{4}{3}-1-\frac{\delta(n)}{2}-\frac{1}{n}+h(n),
\end{array}$$
where $\delta(n)=\left\{\begin{array}{ll}1,&{\rm if}\;n\equiv 0\pmod{2},\\0,&{\rm if}\;n\equiv 1\pmod{2}.\end{array}\right.$

By Lemmas \ref{5002}, \ref{5003}, we have
$${\rm log}R>\varphi(n){\rm log}2-{\rm log}{2n}>c(n){\rm log} 2\cdot n^{\frac{3}{4}}-{\rm log} {2n},$$
$${\rm log}R<{\rm log}4+\delta(n){\rm log}\frac{4}{3}-1-\frac{\delta(n)}{2}-\frac{1}{n}+1.28n^{\frac{1}{4}}.$$
It is easy to calculate that the inequality $${\rm log}4+\delta(n){\rm log}\frac{4}{3}-1-\frac{\delta(n)}{2}-\frac{1}{n}+1.28n^{\frac{1}{4}}>c(n){\rm log} 2\cdot n^{\frac{3}{4}}-{\rm log} {2n}$$ holds for $n\geq 7$ if and only if $$n\in\{7,8,9,10,11,12,13,14,15,16,17,18,19,20,21,22,24,26,30,34,38,42,46,50,54 \}.$$
Hence  the inequality $${\rm log}4+\delta(n){\rm log}\frac{4}{3}-1-\frac{\delta(n)}{2}-\frac{1}{n}+h(n)>\varphi(n){\rm log}2-{\rm log}{2n}$$ holds for $n\geq 7$ if and only if $n\in D=\{8,9,10,12,14,18,20,24,30\}.$   By Lemma \ref{5001}, we can  straightly calculate that there is no $n\in D$ meeting the  assumptions. This completes the proof.   \qed\epf

\vskip 2mm
We are now in the position to prove the main theorem.

\vskip 2mm

{\bf  Proof of Main Theorem }\;\;The sufficiency  is trivial. We need only show the necessity.  We may assume that $p(x)\in\mathbb{F}_q[x]$ is monic and reducible and of degree $n\geq 1.$ Let $$p(x)=p_1(x)^{r_1}\cdots p_k(x)^{r_k}$$ be the standard decomposition, where $p_i(x)$ is monic and irreducible and of degree $e_i\geq 1,$   $r_i\geq 1,$  $1\leq i\leq k.$  By (3) of Proposition  \ref{001}, we have $r_1=r_2=\cdots=r_k=1.$  Hence $$p(x)=p_1(x)\cdots p_k(x),\;\;n=\sum\limits_{i=1}^ke_i,\;\;{\rm and}\;\;\prod\limits_{i=1}^k(q^{e_i}-1)|(q^n-1).$$

 If $q\geq 3,$  then,  by Proposition \ref{002}, we have $q=3,$ $k=2,$ $e_1=e_2=1,$ or  $q=3,$ $k=4,$ $e_1=e_2=e_3=e_4=1.$ But there are only three distinct monic irreducible polynomials of degree one in $\mathbb{F}_3[x],$  hence $p(x)$ is the product of any 2 distinct monic irreducibles of degree $1.$ Hence  
 $$L_{\mathbb{F}_3[x]}=\{ax(x+1), ax(x-1), a(x+1)(x-1)\in\mathbb{F}_3[x], a=1,2\}.$$

 If $q=2,$  then the $e_i's$ satisfy the assumptions of Proposition \ref{5004}, hence we have ${\rm (i)}\; n=2, k=2, e_1=e_2=1;$  or ${\rm (ii)}\; n=4, k=3, e_1=e_2=1,  e_3=2;$ or ${\rm (iii)}\; n=6, k=3, e_1=1, e_2=2, e_3=3.$ On the other hand, the irreducibles of degree 1 are $x$ and $x+1;$ $x^2+x+1$ is the unique irreducible of degree 2; the irreducibles of degree 3 are $x^3+x+1$ and $x^3+x^2+1.$ Hence
$$\begin{array}{ll}L_{\mathbb{F}_2[x]}=&\{x(x+1), x(x+1)(x^2+x+1), x(x^2+x+1)(x^3+x+1),\\
 &\;\; (x+1)(x^2+x+1)(x^3+x+1), x(x^2+x+1)(x^3+x^2+1),\\&\;\; (x+1)(x^2+x+1)(x^3+x^2+1)\in \mathbb{F}_2[x]\}.\end{array}$$ This completes the proof.\qed

\end{document}